\documentclass[12pt, reqno]{amsart}
\usepackage{amsmath, amsthm, amscd, amsfonts, amssymb}

\usepackage{amssymb}
\usepackage{amsfonts}
\usepackage{amsmath}

\usepackage{graphicx}

\usepackage{enumerate}

\usepackage{tikz}
\usetikzlibrary{arrows}

\newtheorem{theorem}{Theorem}[section]

\newtheorem{corollary}[theorem]{Corollary}
\theoremstyle{definition}
\newtheorem{definition}[theorem]{Definition}
\newtheorem{example}[theorem]{Example}

\theoremstyle{remark}
\newtheorem{remark}[theorem]{Remark}
\numberwithin{equation}{section}

\newcommand{\Cstar}{C^\ast}
\newcommand{\C}{$C^\ast$}
\newcommand{\G}{\mathcal{G}}

\usepackage{datetime}

\begin{document}

\title{On graph algebras from interval maps}

\author{Carlos Correia Ramos}
\address{Centro de Investiga\c c\~ao em Matem\'atica e
Aplica\c c\~oes,
 Department of Mathematics,
Universidade de \'{E}vora,
R.\ Rom\~{a}o Ramalho,
59, 7000-671 \'{E}vora,\\
 Portugal}
   \email{ccr@uevora.pt}

\author{Nuno Martins}
\address{Department of Mathematics, CAMGSD,
Instituto Superior T\'{e}cnico, University of Lisbon,
Av.\ Rovisco Pais 1, 1049-001 Lisboa, Portugal}
\email{nmartins@math.tecnico.ulisboa.pt}

\author{Paulo R.\ Pinto}
\address{Department of Mathematics, CAMGSD,
Instituto Superior T\'{e}cnico, University of Lisbon,
Av.\ Rovisco Pais 1, 1049-001 Lisboa, Portugal}
\email{ppinto@math.tecnico.ulisboa.pt}

\begin{abstract}
We produce and study a family of representations of relative graph algebras on Hilbert spaces that
arise from the orbits of points of one dimensional dynamical systems, where the underlying
Markov interval maps $f$ have escape sets. We identify when such representations are faithful
in terms of the transitions to the escape subintervals.
\end{abstract}

\maketitle




\section{Introduction}

For $x$ in the escape set of a Markov interval map $f$, see \cite{RMP10}, we find a relative graph algebra $\Cstar(\G,V_x)$, see  \cite{Tom}, with $\G$ encoding the dynamics of $f$ and a carefully chosen subset $V_x$ of vertex set of $\G$, and a permutative representation $\nu_x$ of this $\Cstar$-algebra on the underlying Hilbert space $H_x$ of the orbit of $x$ such that $\nu_x$ is faithful. This is proved in Theorem \ref{thmnice100}, the main result of this paper. The unitary equivalence within this family of representations is studied through the tree structure of each orbit, see Corollary \ref{lastcor}.

Several concrete subclasses of representations of Cuntz, Cunt-Krieger and graph
algebras have been studied.
We can trace this back to the pioneering work by Bratteli and Jorgensen \cite{BJ},
where they started studying the permutative representations (where the concrete
operators that will generate the representation permute the vectors of the orthonormal
basis of the Hilbert space at hand) of the Cuntz algebras $\mathcal{O}_n$.
From the applications viewpoint, and besides its own right,
applications of representation theory of Cuntz and Cuntz--Krieger
algebras to wavelets, fractals, dynamical systems, see e.g.\
\cite{BJ,marcolli}, and quantum f\/ield theory in~\cite{AK} are
particularly remarkable. For example, it is known that these
representations of the Cuntz algebra serve as a computational tool
for wavelets analysts, see~\cite{Palle2005}, as
such a representation on a Hilbert space $H$ induces a subdivision
of $H$ into orthogonal subspaces. Then the problem in wavelet
theory is to build orthonormal bases in $L^2(\mathbb{R})$ from
these data. Indeed this can be done and these wavelet
bases have advantages over the earlier known basis
constructions (one advantage is the ef\/f\/iciency of computation).
This method has also been applied to the context of fractals that
arise from af\/f\/ine iterated function systems \cite{DuPalle}. Some
of these results have been extended to the more general class of
Cuntz--Krieger algebras, see \cite{CMP, marcolli}, to graph algebras \cite{Daniel} and more recently to higher-rank graph algebras \cite{JorgensenFarsi}.
Much of these applications usually rely on the study of irreducible representations, their
unitary equivalence classes and identify the faithful ones. We remark that if the $\Cstar$-algebra is simple (for example Cuntz algebras) then a nonzero representation is faithful.
So, as soon a new family of permutative representations is build, it is expected that the study of their unitary equivalence of representations will be of great importance.

The relative graph algebras are universal \C-algebras introduced by Muhly and Tomforde in \cite{Tom} as generalizations of graph algebras and Toeplitz algebras. For a finite oriented graph $\G$ and a subset $V$ of vertices of $\G$, the relative graph algebra $\Cstar(\G,V)$ is like the graph algebra $\Cstar(\G)$ except that the relations
$$p_v=\sum_{e: s(e)=v} s_es_e^\ast$$
are required to hold only for $v\in V$. If $V$=$\emptyset$, $\Cstar(\G,\emptyset)$ is the Toeplitz algebra \cite{Rae1} and if $V=\G^0$ is the whole set of vertices, then $\Cstar(\G,\G^0)$ is the graph algebra $\Cstar(\G)$
(which is a Cuntz-Krieger algebra $\mathcal{O}_A$ if $\G$ is the graph represented by an adjacency matrix $A$, see \cite{CK,Rae2}).

In this paper we investigate how dynamical systems arising from interval
maps can produce relative graph algebras. This is done by yielding representations of relative graph algebras on the orbits of the underlying dynamical system. Moreover, we study when such representations are faithful. This provides a situation in which relative graph algebras prove to be convenient. We follow the spirit of \cite{CMP,RMP5,RMP10} where this was pursued for the two extremal cases, the Cuntz-Krieger algebra $\Cstar(\G,\G^0)$ and the Toeplitz algebra $\Cstar(\G,\emptyset)$ case. Since $\Cstar(\G,V)$ is in general nonsimple, faithfulness of the representations do not come for free.

Indeed, the underlying interval maps $f: \hbox{dom}(I)\to I$ (with dom$(f)\subset I$) can have escape sets, where the orbits of the points may fall outside of the domain dom$(f)$ of $f$. This originates an {\it escape transition matrix} $\widehat{A}_f$  as in \cite{RMP10,RMP11} that encodes not only the transitions from Markov subintervals to Markov subintervals (which are recorded in a matrix $A_f$) but also the transitions to escape subintervals of $I$.
If we gather the Markov subintervals first and then the escape subintervals, then there exists a permutation matrix $P$ such that
\begin{equation}
P\widehat{A}_fP^T= \left[
\begin{array}{c|c}
A_f & B_f \\
\hline
0_{m\times n} & 0_{m\times m}
\end{array}
\right]
\label{PAPT}
\end{equation}
where the lower blocks are matrices with all entries equal to zero and the matrix $B_f$ collects the transitions from Markov subintervals into escape subintervals (see \cite{RMP11}).

More precisely, in \cite{RMP10} we considered the interval maps $f$ where the orbit of a point
$x$ can be in the \textit{escape set} of $f$ -- see \cite{Fal} -- namely, $f^{k}(x)\in I$
does not belong to the domain of $f$ for a certain $k\in \mathbb{N}$. We
likewise define a Hilbert space $H_{x}$ from the backward orbit of $x$ and
partial isometries on $H_{x}$. In that situation we proved
that these operators do define a representation $\nu_{x}$ of the Toeplitz algebra associated to the
transition matrix $A_f$ on $H_{x}$ (which is no longer a
representation of the Cuntz-Krieger algebra $\mathcal{O}_{A_{f}}$ if the above matrix $B_f\not=0$).
This representation $\nu_x$ is in general non faithful and this is where the relevance of the relative graph algebras will play a role as we see in the sequel.



Given such an interval map $f$, we concentrate on the graph $\G_{A_f}$ attached to the matrix $A_f$ and work on the family of relative graph algebras we can associate to $\G_{A_f}$. If the point $x$ is in the escape set of $f$, then under iteration $x$ will fall into one of the escape sets of $f$, which is encoded in some column $j$ of $B_f$. Then we consider the set $V_x$ of rows of $A_f$ for which there is no transition to the escape set $E_j$. Then we prove that indeed we can define a representation $\nu_x$ of the relative graph algebra $\Cstar(\G_A,V_x)$ on the Hilbert space $H_x$ and show that $\nu_x$ is faithful. However $\nu_x$ is also a representation of the Toeplitz algebra $\Cstar(\G_A,\emptyset)$ or even a representation of any $\Cstar(\G_{A_f},V)$ with $V\subset V_x$, so we may apply \cite{RMP10} which leads us to the study of the irreducibility and unitary equivalence of this family $\{\nu_x\}_{x\in I}$ of representations.

The plan for the rest of the paper is as follows. In Sect.\ \ref{secmarkovescape} we first review the dynamical systems background for the interval maps with escape sets, emphasizing the notion of {\it escape transition matrix} as in Definition \ref{defhatA}.

In Sect.\ \ref{secgraphalgs}, we first review some elements of the (relative) graph algebra theory, and then proceed in Subsect.\ \ref{subrepsescape} with a construction of a representation of relative graph algebras for every point $x$ and interval map $f$. The underlying Hilbert space $H_x$ encodes the (generalised) orbit of $x$, and the generating partial isometries of the graph algebra are defined as in \eqref{eqsijpi}.
If $A_f$ is the transition matrix of an interval map $f$,  $\G_{A_f}$ its oriented graph, then fixing $V\subseteq  \G_{A_f}^0$ our main result, Theorem \ref{thmnice100}, shows how we obtain a representation $\nu_x$ of the \C-algebra $\Cstar(\G_{A_f},V)$ on $H_x$ and the conditions under which $\nu_x$ is faithful.
Example \ref{example111} guides us through a concrete example where we build an interval map (and we also write the transitions matrices $A_f$ and $\widehat{A}_f$ together with the graph $\G_{A_f}$) -- then Theorem \ref{thmnice100} tells us when $\nu_x$ is a representation and when $\nu_x$ is a faithful representation of $\Cstar(\G_{A_f},V)$, for every set $V\subset \G_{A_f}^0$
We show in Corollary \ref{corfaithreps} the richness of this class of representations.
We conclude the paper with some remarks that are consequences of the above results together with \cite{RMP10}.

\section{Markov interval maps with escape sets}
\label{secmarkovescape}

Given $n\in \mathbb{N}$, let
$\Gamma$=$\{c_{0},c_{1}^{-},c_{1}^{+},...,c_{n-1}^{-},c_{n-1}^{+},c_{n}\}$
be an ordered set of (at most) $2n$ real numbers such that
\begin{equation}
c_{0}<c_{1}^{-}\leq c_{1}^{+}<c_{2}^{-}\leq ...<c_{n-1}^{-}\leq
c_{n-1}^{+}<c_{n}\text{.}
\label{eq2}
\end{equation}

Given $\Gamma $ as above, we define the collection of closed intervals
$C_{\Gamma }=\{I_{1},...,I_{n}\}$, with
\begin{equation}
I_{1}=\left[ c_{0},c_{1}^{-}\right],\ ...,\ I_{j}=\left[
c_{j-1}^{+},c_{j}^{-}\right],..., I_{n}=\left[c_{n-1}^{+}, c_{n}\right].
\label{eq3}
\end{equation}%
We also consider the collection of open intervals $\{E_{1},...,E_{n-1}\}$,
with
\begin{equation}
E_{1}=\left] c_{1}^{-},c_{1}^{+}\right[,\ ...,\ E_{n-1}=\left]
c_{n-1}^{-},c_{n-1}^{+}\right[,  \label{eq4}
\end{equation}%
in such a way that $I:=[c_{0},c_{n}]=\left(\cup _{j=1}^{n}I_{j}\right)
\bigcup \left(\cup _{j=1}^{n-1}E_{j}\right)$.


We now consider the interval maps for which we can construct partitions of
the interval $I$ as in \eqref{eq2}, \eqref{eq3} and \eqref{eq4}.

\begin{definition}[See \cite{CMP}]
Let $I\subset \mathbb{R}$ be an interval. A map $f$ is  called a Markov interval map or it is
in the class $\mathcal{M}(I)$ if it satisfies the following properties:

\begin{enumerate}
\item[(P1)] \textrm{[}Existence of a finite partition in the domain of $f$%
\textrm{]} There is a partition $C=\left\{I_{1},...,I_{n}\right\}$ of closed
intervals with $\#\left(I_{i}\cap I_{j}\right)\leq 1$ for $i\neq j$,
$\mathrm{dom}(f)=\bigcup _{j=1}^{n}I_{j}\subset I$, $\min(I_1)=\min(I)$, $\max(I_n)=\max(I)$  and $\mathrm{im}(f)=I$.

\item[(P2)] \textrm{[}Markov property\textrm{]} For every $i=1,...,n$ the
set $f(I_{i})\cap \left(\bigcup_{j=1}^{n}I_{j}\right)$ is a non-empty union
of intervals from $C$.


\item[(P3)] \textrm{[}Expansive map\textrm{]} $f_{|I_{j}}\in \mathcal{C}%
^{1}(I_{j})$, monotone and $|f_{|I_{j}}^{\prime }(x)|>b>1$, for every $x\in
I_{j},j=1, ..., n$, and some $b$.

\item[(P4)] \textrm{[}Aperiodicity\textrm{]} For every interval $I_{j}$ with
$j=1,...,n$ there is a natural number $q$ such that $\mathrm{dom}(f)\subset
f^{q}(I_{j})$.
\end{enumerate}
\label{def1}
\end{definition}

The minimal partition $C$ satisfying the Definition \ref{def1} is denoted by $C_{f}$. We remark that
the Markov property (P2) allows us to encode the transitions between the
intervals in the so-called (Markov) transition $n\times n$ matrix
$A_{f}=(a_{ij})$, defined as follows:
\begin{equation}
a_{ij}=\left\{
\begin{array}{l}
1\text{ if }f(\mathring{I}_{i})\supset \mathring{I}_{j}, \\
0\text{ otherwise}
\end{array}
\right.  \label{mx}
\end{equation}
where $\mathring{A}$ denotes the interior of the set $A$.
A map $f\in \mathcal{M}(I)$ uniquely determines (together with the minimal
partition $C_{f}=\left\{ I_{1},...,I_{n}\right\}$):

\begin{enumerate}[(i)]
\item The $f$-invariant set $\Omega _{f}:=\{x\in I:\ f^{k}(x)\in \mathrm{dom}%
(f)$ for all $k=0,1,...\}$.

\item The collection of open intervals $\left\{ E_{1},...,E_{n-1}\right\} $,
such that ${I\setminus \bigcup _{j=1}^{n}I_{j}}=\bigcup _{j=1}^{n-1}E_{j}$.

\item The transition matrix $A_{f}=\left( a_{ij}\right)_{i,j=1,..,n}$.
\end{enumerate}

Let $I$ be an interval and $f\in \mathcal{M}(I)$. Once we chose a
Markov partition $C$, then there exists a unique set $\Gamma$ of boundary points
satisfying \eqref{eq2}, \eqref{eq3} and \eqref{eq4}. Note that the Markov property implies that
$f\left(\Gamma \right) \subset \Gamma$.

Matrices $A$ for which there exists a positive integer $m$ such that all the entries of $A^m$
are non-zero are called primitive.
We note that the matrix $A_f$ is primitive (thus irreducible) because $f\in
\mathcal{M}(I)$, see Definition \ref{def1}.

Note that we explicitly use the notation $c_{j}^{+}$ and $c_{j}^{-}$ to
represent the boundary points of the Markov partition of $f$, whose domain equals the union of the closed subintervals $I_{j}$. This is done to
formalize important notions for the study of the orbit structure of $f$, see \cite{Fal}. If
 $c_{j}^{-}<c_{j}^{+}$ then the points $c_{j}^{-}$ and $c_{j}^{+}$ are distinct. In this case, the domain
of $f$ is disconnected and the escape interval $E_{j}=\left]
c_{j}^{-},c_{j}^{+}\right[$ is non-empty.
In the case $c_{j}^{-}$ and $c_{j}^{+}$ are equal, then $c_{j}^{+}$, $c_{j}^{-}$ represent the side
limits of the point $c_{j}$, and the escape interval $E_{j}$ is empty.
In this case, $f$ can be either
differentiable, continuous but non-differentiable or discontinuous at $c_{j}$, see \cite{Bandeira}. In
this last case, we could arbitrarily choose the side limit $c_{j}^{-}$ or
$c_{j}^{+}$ where the map $f$ is defined, since $f\left( c_{j}^{-}\right)
\neq f\left( c_{j}^{+}\right)$. However the two side limits give origin to
important orbits (the orbit of $f(c_j^-)$ and $f(c_j^+)$) for the characterization of the dynamics.
In the present paper these orbits of $f\left( c_{j}^{-}\right)$ and $f\left( c_{j}^{+}\right)$ will not be
used as we will see in the sequel.

Note that $\Omega _{f}$ is the set of points that remain in $\mathrm{dom}(f)$
under iteration of $f$, and is usually called a \textit{cookie-cutter set},
see \cite{Fal}. The open set
\begin{equation}
E_{f}:=I\setminus
\Omega_{f}=\bigcup_{k=0}^{\infty }f^{-k}\left(\bigcup _{j=1}^{n-1}E_{j}\right)
\end{equation}
is called the \textit{escape set}. Every point in $E_f$ will eventually
fall, under iteration of $f$, into some interval $E_{j}$
(where $f$ is not defined) and the iteration process ends.
We may say that $x$ is in the escape set $E_{f}$ of $f$ if and only if there is $k\in
\mathbb{N}$ such that $f^{k}\left(x\right) \notin \mathrm{dom}(f)$.

We define an equivalence relation $R_f$ on $I$ by
\[R_f:=\{(x,y): f^n(x)=f^m(y) \text{ for some } n,m \in \mathbb{N}_0 \}. \]
(The restriction of this equivalence relation to $\Omega_f$ was considered in \cite{CMP}.)
The relation
$R_{f}$ is a countable equivalence relation in the sense that the equivalence
class $R_{f}(x)$ of $x\in I$, is a countable set, and
can be seen as the generalized orbit of $x$ (forward and backward orbit).

Given a set $X\subset \hbox{dom}(f)$, $R_{f}\left( X\right) $ is the set of
generalized orbits of the points in $X$.
We are interested in studying points in the escape
set or in the so-called {\it regular set}
\begin{equation}
\Lambda_{f}=\Omega _{f}\setminus R_{f}\left(\Gamma
\right)
\label{eqregular}
\end{equation}
of $f\in\mathcal{M}(I)$.  The regular set $\Lambda_f$ is the set of the points in the maximal invariant set, under $f$, which are not boundary points and are not pre-images of boundary points. Note that
since we are in Markov case the image of a boundary point is always a
boundary point.

Let $f\in \mathcal{M}(I)$ be such that $E_{f}\neq \emptyset$. This means
that there is at least one non empty open interval
$E_{j}=\left]c_{j}^{-},c_{j}^{+}\right[ $,
with $c_{j}^{-}\neq c_{j}^{+}$, with $j\in\{1,...,n-1\}$.
The non-empty open subinterval $E_{j}$ is called an \textit{escape interval}.

In order to describe symbolically the escape orbits, we extend the symbol
space adding a symbol for each escape interval $E_{j}$, which will represent
an end for the symbolic sequence. For each escape interval $E_{j}$ we
associate a symbol $\widehat{j}$ to distinguish from the symbol associated to
the interval $I_{j}$ of the partition. That is, we consider the symbols ordered by:
\begin{equation}
1<\widehat{1}<2<\widehat{2}< ... < n-1< \widehat{n-1} <n\text{.}
\label{eqorder}
\end{equation}

If $E_{j}$ is not an interval, that is $E_{j}=\emptyset $, then there is no symbol
$\widehat{j}$. Moreover, we define
\begin{equation}
\Sigma_{E_f}=\left\{\widehat{j}:  E_{j}\neq \emptyset,\ j\in\left\{1,..., n-1\right\}
\right\}.
\label{eqhatsymbols}
\end{equation}

For
every $y\in E_{f}$ there is a least natural number $\tau \left( y\right) $
such that $f^{\tau \left(y\right)}\left(y\right)\notin\mathrm{dom}\left(
f\right)$, which means that, $f^{\tau \left( y\right)}\left(y\right)\in
E_{j}$, for some $j$ such that $E_{j}\neq\emptyset$.
The final escape point, for the orbit of $y$, is then denoted by
$e\left( y\right) :=f^{\tau \left( y\right)}\left( y\right) $ and the final
escape interval index is denoted by $\iota \left(y\right) $, that is, if
$f^{\tau \left( y\right) }\left( y\right) \in E_{j}$ then $\iota \left(
y\right) =\widehat{j}$.

Thus we have an index set $\{1,...,n\}\bigcup\Sigma_{E_f}$ which is ordered as in \eqref{eqorder} and \eqref{eqhatsymbols}.
In order to deal with the possible transitions from Markov transition intervals to
escape intervals we define the escape transition matrix $\widehat{A}_{f}$ as follows.

\begin{definition}[see \cite{RMP10, RMP11}]
Given the transition matrix $A_f$ as in \eqref{mx}, we define a matrix
 $\widehat{A}_{f}=\left( \widehat{a}_{ij}\right)$ indexed by $\{1,...,n\}\cup\Sigma_{E_f}$ such that
\begin{equation}
\widehat{a}_{ij}=\left\{
\begin{array}{ll}
a_{ij}&\text{ if } i,j\in\{1,...,n\}, \\
1&\text{ if } i\in \{1,...,n\},\ j=\hat{k}\in \Sigma _{E_f}\ \hbox{and}\ \mathring{I}_{i}\cap f^{-1}\left(
E_{k}\right) \neq \emptyset,\\
0&\text{  }  \hbox{otherwise}.
\end{array}%
\right.
\label{mmx}
\end{equation}
\label{defhatA}
\end{definition}
For row and column labeling, the matrix $\widehat{A}_f$ is defined by considering the order in \eqref{eqorder}.
Note that if we use the row and column labeling order $I_1...I_n E_1...E_m$ with $\#\Sigma_{E_f}=m$, then
we obtain the reordered transition matrix in the form
\begin{equation}
\left[
\begin{array}{c|c}
A_f & B_f \\
\hline
0_{m\times n} & 0_{m\times m}
\end{array}
\right],
\label{eqmatrixa00}
\end{equation}
where $A_f$ is the Markov transition matrix of $f$ and $B_f$ is the transition matrix from Markov subintervals to the escape subintervals. We write $0_{p\times q}$ for the $p\times q$ matrix with zeros everywhere, whereas $1_{p\times q}$ is the $p\times q$ matrix with ones everywhere.

\section{Relative graph $C^\ast$-algebras from interval maps}
\label{secgraphalgs}

\subsection{Background on relative graph algebras}

Let $\G=(\G^0,\G^1,r,s)$ be a (oriented, finite) graph, with $\G^0$ the vertex set, $E^1$ the edge set and $r,s: \G^1\to \G^0$ the range and source maps, respectively. The graph algebra $\Cstar(\G)$ is the universal \C-algebra generated by partial isometries $\{s_e: e\in \G^1\}$ with commuting range projections together with a collection of mutually orthogonal projections $\{p_v: v\in \G^0\}$ that satisfy the relations
$$s_e^\ast s_e= p_{r(e)},\quad p_v=\sum_{e: s(v)=v} s_e s_e^\ast\ \hbox{whenever}\ s^{-1}(v)\not=\emptyset.$$
We remak that then have $s_e s_e^\ast \leq p_{s(e)}$ since we have assumed $\G$ is a finite graph.
We refer to the textbook \cite{Rae} for the standard theory of graph algebras, but \cite{Rae} uses a different convention to the one you are using, namely that the partial isometries representing the edges go in the same direction as the edges (cf.\ \cite[p.1 2 Conventions]{Rae}).

For any
$n\times n$ matrix $A=(a_{ij})$ with entries in $\{0,1\}$, we can construct a
directed graph $\G_A=(\G_A^1,\G_A^0, r,s)$ with
\begin{equation}
\G_A^0=\{1,...,n\},\ \G_A^1=\{e_{ij}: i,j\in\G_A^0, a_{ij}= 1\}\quad \hbox{with}\ s(e_{ij})=i,\ r(e_{ij})=j
\label{eqGA}
\end{equation}
(i.e.\ we draw an edge $e_{ij}$ from $i$ to $j$ if and only if $a_{ij}=1$)
where $s$ and $r$ are the source and range maps, respectively. The Cuntz-Krieger algebra $\mathcal{O}_A$ is then the graph \C-algebra $\Cstar(\G_A)$, see \cite[Prop.\ 4.1]{Rae2}. The particular case, but important, when $A$ is full ($a_{ij}=1$ for all $i,j$), the Cuntz-Krieger algebra $\mathcal{O}_A$ is the Cuntz algebra $\mathcal{O}_n$, where the partial isometries are all isometries.

For any $V\subset \G^0$, the {\it relative graph algebra} $\Cstar(\G,V)$ (see \cite[Def.\ 3.4]{Tom}) is the universal \C-algebra generated by partial isometries $\{s_e: e\in \G^1\}$ and mutually orthogonal projections $\{p_v: v\in \G^0\}$ satisfying the following set of relations $\mathcal{R}_V$:

\begin{eqnarray}
  s_e^\ast s_e &=&  p_{r(e)}, \label{eqsesepr}\\
   s_e s_e^\ast &\leq& p_{s(e)},\ \hbox{for all}\ e\in \G^1, \ \hbox{and} \label{eqseseps}\\
  p_v &=& \sum_{e\in \G^1:\ s(e)=v} s_e s_e^\ast, \label{eqpv=}
\end{eqnarray}
 for all $v\in V$ such that $\{e\in \G^1: s(e)=v\}\not=\emptyset$.
As remarked in \cite{Tom}, for $v\in \G^0\setminus V$ we have $p_v> \sum_{s(e)=v} s_e s_e^\ast$. Besides this, $C^\ast(\G,V_1)$ is a quotient of $C^\ast(\G,V_2)$ whenever $V_2\subset V_1$, and in particular $C^\ast(\G,V)$ is a quotient of $C^\ast(\G,\emptyset)$, for any $V$. Namely, if $J$ is the ideal generated by the {\it gap} projections $$\{p_v-\sum_{e\in \G^1:\ s(e)=v} s_e s_e^\ast: v\in V_1\}$$
then $C^\ast(\G,V_1)=C^\ast(\G,V_2)/J$.

The relative graph algebra $\Cstar(\G,\emptyset)$ is the so-called Toeplitz algebra \cite{Rae1,Tom}, and denoted by $\mathcal{T}_A$ in \cite{RMP10}, where $A$ is a 0-1 matrix whose associated graph is $\G_A$, see \eqref{eqGA}.

The Cuntz-Krieger algebra $\mathcal{O}_A$ of a primitive matrix $A$ is known to be simple (thus the Cuntz algebra $\mathcal{O}_n$ is simple), but the Toeplitz algebra $\mathcal{T}_A$ is non-simple, see
\cite[Corollary 4.3]{Rae1}.
Therefore every non-zero representation of $\mathcal{O}_A$ is automatically faithful, whereas for the Toeplitz algebra or other relative graphs algebras $\Cstar(\G,V)$, faithfulness does not come for free.

\subsection{Representations of relative graph algebras arising from interval maps}
\label{subrepsescape}

Now assume that an $(n+m)\times (n+m)$ matrix $\widehat{A}_f$ with entries in $\{0,1\}$ is given produced from a Markov interval map $f\in \mathcal{M(I)}$ whose transition matrix is an $n\times n$ primitive transition matrix $A_f$. Thus
$$P\widehat{A}_f P^T=\left[
\begin{array}{c|c}
A_f & B_f \\
\hline
0_{m\times n} & 0_{m\times m}
\end{array}
\right],$$
with $P$ a permutation matrix, $B_f$ an $n\times m$ matrix that encodes the transitions from the Markov intervals into the escape intervals.

Let $H_{x}=\ell^2(R_f(x))$, with $x\in E_{f}$, be the Hilbert space with canonical base
\begin{equation*}
\left\{ \left\vert z\right\rangle: f^{k}\left(z\right)= e\left(x\right)\ \hbox{for some}\ k\in \mathbb{N}_{0}\right\}.
\end{equation*}
 Note that there is a
special vector basis which is $\left\vert e\left(x\right)\right\rangle $.
The rank one projection on the 1-dimensional space $\mathbb{C}\left\vert z\right\rangle$
is denoted by $P_{z}$, or as usual in Dirac notation, $P_{z}=\left\vert
z\right\rangle \left\langle z\right\vert $.

Let $x\in E_{f}$ and let $T_{i},$ $i=1,...,n$ be defined by:
\begin{equation}
T_{i}\left\vert y\right\rangle =\chi _{f(I_{i})}(y)\left\vert
f_{i}^{-1}(y)\right\rangle \text{ for }y\in R_{f}\left(x\right)
\label{defTii}
\end{equation}%
where $f_i:=f|_{I_i}$, $\chi_{B}$ denotes the characteristic function on a set $B$. The following equations are naturally
satisf\/ied, and will be used in the sequel,
\begin{equation}
\left\vert f\circ f_{i}^{-1}\left( y\right) \right\rangle =\chi _{f\left(
I_{i}\right) }\left( y\right) \left\vert y\right\rangle, \ \ \ \left\vert
f_{i}^{-1}\circ f_{|I_{i}}\left( y\right) \right\rangle =\chi _{I_{i}}\left(
y\right) \left\vert y\right\rangle
\label{eqnice11}
\end{equation}
where $f_i$ is invertible $f|_{I_i}^{-1}: f(I_i)\to I_i$ because $f$ is a Markov map (note that $f\circ f_{i}^{-1}\left(y\right) $ is not defined for $y\notin f\left( I_{i}\right)$. However we set
$\left\vert f\circ f_{i}^{-1}\left( y\right) \right\rangle =0$ so that
the above expression is well defined.)

Its adjoint, $T_{i}^{\ast }$, is given by
\begin{equation*}
T_{i}^{\ast }\left\vert y\right\rangle =\chi _{I_{i}}(y)\left\vert
f(y)\right\rangle \text{.}
\end{equation*}
In particular, $T_{i}^{\ast }\left\vert e(x)\right\rangle=0$ for all $i=1,...,n$,
and $T_{i}\left\vert e(x)\right\rangle =0$ if there are no transitions from
the interval $I_{i}$ to the escape interval $E_{j}$. If there is any
transition, a $z\in I_i$ with $f(z)=e(x)$ say,  then $T_{i}\left\vert e(x)\right\rangle =\left\vert z\right\rangle$
such that $z\in I_{i}$ and $f\left( z\right)=x$.
The operators $T_1,...,T_n$, are all partial isometries.

In order to translate these partial isometries into the graph algebras framework,
we take the associated graph $\G_{A_f}$ and define the following operators:
\begin{equation}
S_{e_{ij}}:= T_iT_jT_j^\ast\ \hbox{if}\ e_{ij}\in\G_{A_f}^1,\ \hbox{and}\ P_i:=T_i T_i^\ast,\ Q_i:=T_i^\ast T_i\quad  i,j=1,...,n.
\label{eqsijpi}
\end{equation}
We will write $S_{ij}$ instead of $S_{e_{ij}}$. If we use the definition of the partial isometries $T_1,...,T_n$, for vectors $|y\rangle$ such that $y\in R_f(x)$, we obtain
\begin{eqnarray}
  S_{ij} |y\rangle &=& \chi_{I_j}(y)\ | f_i^{-1}(y)\rangle,\quad (\hbox{if}\ a_{ij}=1), \label{eqSij}\\
  P_i|y\rangle&=&  \chi_{I_i}(y)\ |y\rangle, \label{eqPi}\\
  Q_i|y\rangle &=& \chi_{f(I_i)}(y)\ |y\rangle.
\end{eqnarray}
It is clear that $S_{ij}$ are nonzero partial isometries for $e_{ij}\in\G_{A_f}^1$ (i.e.\ $a_{ij}=1$) with
$$S_{ij}^\ast |y\rangle =\chi_{I_j}(f(y))\chi_{I_i}(y)\ |f (y)\rangle,$$
and $P_i$ nonzero projection for all $i\in \G_{A_f}^0$ (i.e.\ $i=1,...,n$). We also note that the range projections $P_1,..., P_n$ are pairwise orthogonal because $ \mathring{I}_i\cap \mathring{I}_j=\emptyset$ for $i\not= j$ and $Q_iP_j=a_{ij} P_j$.

Recall from \eqref{eqregular} that $\Lambda_f$ is the regular set of a Markov map $f$. We then remark that if $x\in E_f\cup \Lambda_f$, then $\Gamma\cap R_f(x)=\emptyset$.
\begin{theorem}
Let $A_f$ be the transition matrix of an interval map $f$ and $\widehat{A}_f$ its escape transition matrix so that $E_f\not=\emptyset$.
 Consider its oriented graph $\G_{A_f}$ and fix $V\subseteq  \G_{A_f}^0$. Let $x\in I$.
 \begin{enumerate}
 \item Let $x\in E_f$.
 \begin{enumerate}
 \item If $\widehat{a}_{i \iota(x)}=0$ for $i\in V$, then $s_{e_{ij}}\mapsto S_{ij}$ and $p_i\mapsto P_i$ defined in \eqref{eqsijpi} yield a representation, $\nu_x$, of $\Cstar(\G_{A_f}, V)$ on the Hilbert space $H_x$.
 \item If $\widehat{a}_{i \iota(x)}=0$ for $i\in V$ and $\widehat{a}_{k\, \iota(x)}=1$ for $k\notin V$, then the representation in part (1) is faithful.
\end{enumerate}

\item If $x\in\Lambda_f$, then the operators defined in \eqref{eqsijpi} yield a faithful representation of
 $\Cstar(\G_{A_f}, \G_{A_f}^0)$ on the Hilbert space $H_x$, which is the Cuntz-Krieger algebra $\mathcal{O}_{A_f}$.
\end{enumerate}
\label{thmnice100}
\end{theorem}
\begin{proof}
To prove part (1a), we use the universality of $\Cstar(\G_{A_f},V)$ and check that the nonzero projections $P_i$ and partial isometries $S_{ij}$ do satisfy the defining relations \eqref{eqsesepr}, \eqref{eqseseps} and \eqref{eqpv=}.
First note that $S_{ij}$ is defined only if $a_{ij}=1$ and in this case $S_{ij}=T_iP_j$, so
$$S_{ij}^\ast S_{ij}=P_j T_i^\ast T_i P_j=P_jQ_iP_j=a_{ij}P_jP_j=P_j.$$
We therefore obtain the relation \eqref{eqsesepr}.

It is straightforward to check that $P_i (S_{ij}S_{ij}^\ast)=S_{ij}S_{ij}^\ast$, thus $S_{ij}S_{ij}^\ast\leq P_i$, which implies \eqref{eqseseps}.

Next, we compute
\begin{eqnarray*}
\sum_{e_{ij}\in \G_A^1: s(e_{ij})=i} \ S_{ij}S_{ij}^\ast&=&\sum_{j: a_{ij}=1} T_i T_j T_j^\ast T_i^\ast  =
T_i\left(\sum_{j=1}^n a_{ij} T_jT_j^\ast\right)T_i^\ast\\
 &=&T_i(T_i^\ast T_i-\widehat{a}_{i \iota(x)}P_{e(x)})T_i^\ast
\end{eqnarray*}
where we use \cite[Lemma 4.1]{RMP10} in the last equality. So if $\widehat{a}_{i \iota(x)}=0$, then the RHS of the above expression can be further simplified as follows:
$$T_i(T_i^\ast T_i-\widehat{a}_{i \iota(x)}P_{e(x)})T_i^\ast=T_iT_i^\ast T_iT_i^\ast=P_i=P_{s(e_{ij})}$$
which implies the equality \eqref{eqpv=}.

Now we prove statement (1b).
For such data $\G_A$ and $V$, let us consider the extended (oriented) graph ${\G}_A^V$ taking $\G_A$ and adding a sink $v^\prime$ for each $v\in \G_A^0\setminus V$ as well as edges to this sink from each vertex that feeds (in $\G_A$) into $v$.
Thus, the condition (L) -- every cycle in the graph has an exit -- of ${\G}_A^V$ is inherited from $\G_A$. (Condition (L) holds for $\G_A$ because $A$ is a primitive matrix.)

Let $\{p_v: v\in \G_A^0\}\cup \{p_e: e\in \G_A^1\}$ be the generators of $\Cstar(\G_A,V)$ and $\{\overline{p}_v: v\in ({\G}_A^V)^0\}\cup \{\overline{p}_e: e\in ({\G}_A^V)^1\}$ the generators of the graph algebra $\Cstar({\G}_A^V)$.
Then \cite[Theorem 3.7]{Tom} shows that there is a canonical isomorphism $\phi: \Cstar({\G}_A^V)\to \Cstar(\G,V)$ such that
$$\phi(\overline{p}_v)=\left\{
\begin{array}{lll}
p_v&\hbox{if}& v\in  V\\
&\cr
\sum_{e\in \G_A^1: s(e)=v} s_es_e^\ast &\hbox{if}& v\in \G_A^0\setminus V\\
&\cr
p_v-\sum_{e\in \G_A^1: s(e)=v} s_es_e^\ast &\hbox{if}& v=u^\prime\ \hbox{for some}\ u\in \G_A^0\setminus V.
\end{array} \right.$$
Of course $\nu_x$ is faithful if and only if $\nu_x\circ \phi$ is faithful. By \cite[Theorem 3.1]{Bates},
$\nu_x\circ \phi$ is faithful if the images of the projections $\{\overline{p}_v: v\in ({\G}_A^V)^0\}$ of  $\Cstar({\G}_A^V)$ in $B(H_x)$ are all nonzero. Clearly $P_i\not=0$ for $i\in V$. If $i\not\in V$,
\begin{equation}\label{eqinjnn}
P_i-\sum_{e \in \mathcal{G}_A^1 : s(e) = i} S_{ij}S_{ij}^* = T_i(T^*_iT_i-\sum_{j=1}^n a_{ij}T_jT_j^*)T_i^* = \hat{a}_{i\iota(x)} P_{e(x)} = P_{e(x)}
\end{equation}
where we use \cite[Lemma 4.1]{RMP10} in the last but one equality and the hypothesis $\hat{a}_{i\iota(x)}=1$ for $i \notin V$. Thus $P_i-\sum_{e \in \mathcal{G}_A^1 : s(e) = i} S_{ij}S_{ij}^* \neq 0$ for $i \notin V$. Finally, if $\sum_{e \in \mathcal{G}_A^1 : s(e) = i} S_{ij}S_{ij}^* = 0$ for
$i \notin V$ then Eq. (3.11) would implies that $P_i=P_{e(x)}$ which is not true.
Therefore $\sum_{e \in \mathcal{G}_A^1 : s(e) = i} S_{ij}S_{ij}^* \neq 0$ for $i \notin V$. This finishes the proof of part (1b).

For part (2), we either follow  part (1) or use \cite[Thm.\ 6]{CMP} where it is shown that the partial isometries $T_1,..., T_n$ as in  Eq.\ (\ref{defTii}) produce a faithful representation of the Cuntz algebra $\mathcal{O}_{A_f}$.
\end{proof}

\begin{definition}
We denote by $\nu_x$ (or $\nu_{x,V}$) the representation of
the \C-algebra $\Cstar(\G_{A_f},V)$ produced in Theorem \ref{thmnice100}.
\label{defnux}
\end{definition}

\unitlength 0.1mm
\begin{figure}[ht]
\centerline{\includegraphics[scale=.480]{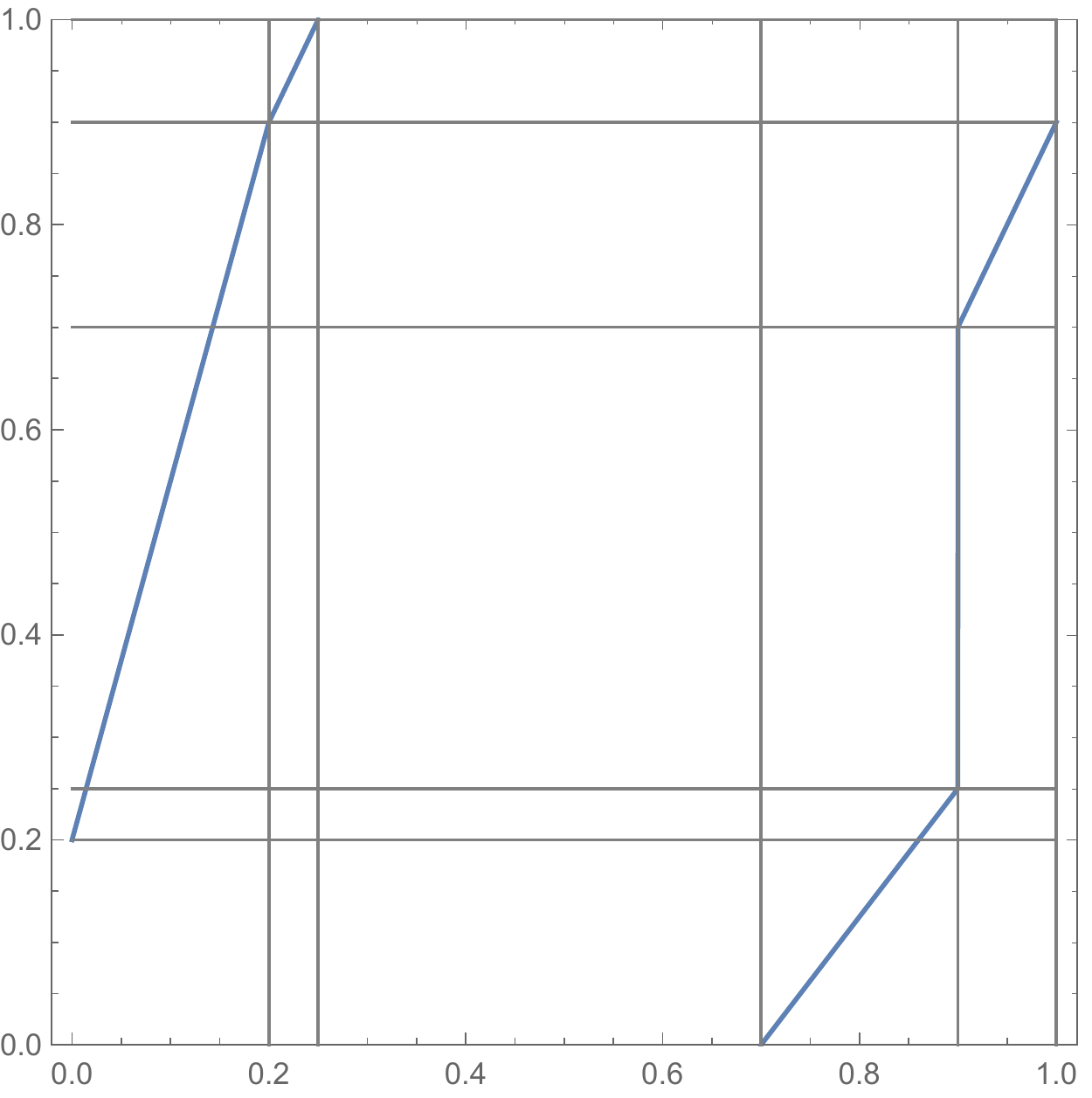}}
 \caption{Graph of the function $f$ in (\ref{eqexample})}
\label{fig1}%
\end{figure}

\begin{example}
Let $f$ be an interval map such that its transition matrix is
$A_f=\left[
\begin{array}{cccc}
0 & 1 &1 & 0 \\
0 & 0 & 0 & 1 \\
1 & 1 & 0 & 0 \\
0 & 0 & 1 & 0
\end{array}%
\right].$ The alphabet of the dynamics has four symbols $\{1,2,3,4\}$.
The escape matrix can be chosen as
$$\widehat{A}_f=\left[
\begin{array}{ccccc}
0 & 1 &1 & 1& 0 \\
0 & 0 & 0 &0& 1 \\
0 & 0 & 0 & 0&0\\
1 & 1 & 0 & 0 &0\\
0 & 0 & 1 & 1&0
\end{array}%
\right].$$
The existence follows from \cite{RMP11} or by constructing a interval map with the requested transitions.
For example,
let $f$ be the interval map defined by (see Fig.\ \ref{fig1})
\begin{equation}\label{eqexample}
f\left(x\right)=\left\{
\begin{array}{c}
\frac{7}{2}x+\frac{1}{5}\text{ \ \ \ \ \ if \ \ \ }x\in I_{1} \\
2x+\frac{1}{2}\text{ \ \ \ \ \ \ if \ \ \ }x\in I_{2} \\
\frac{5}{4}x-\frac{7}{8}\text{ \ \ \ \ \ if \ \ \ }x\in I_{3} \\
2x-\frac{11}{10}\text{ \ \ \ \ \ if \ \ \ }x\in I_{4}
\end{array}
\right.
\end{equation}
with $I_{1}=\left[0,\frac{1}{5}\right]$, $I_{2}=\left[\frac{1}{5},\frac{1}{4}\right]$,
$E_{2}=\left] \frac{1}{4},\frac{7}{10}\right[$, $I_{3}=\left[\frac{7}{10},\frac{9}{10}\right] $, and
$I_{4}=\left[ \frac{9}{10},1\right]$.
Note that
$$\Gamma =\{ 0,\frac{1}{5},\frac{1}{4},\frac{7}{10},\frac{9}{10},1\},\ \ \mathrm{dom}(f) =\bigcup _{i=1}^{4}I_{i}\ \ \hbox{and}\ \mathrm{%
im}\left( f\right) =\left[ 0,1\right] =\bigcup _{i=1}^{4}I_{i}\cup E_2,$$
$f(\frac{9}{10}^-)=\frac{1}{4}$ and $f(\frac{9}{10}^+)=\frac{7}{10}.$
It is now easy to check that  $f\in \mathcal{M}([0,1])$,
$$f(I_1)=I_2\cup E_2\cup I_3,\ f(I_2)=I_4,\ f(I_3)=I_1\cup I_2\ \hbox{and}\ f(I_4)=E_2\cup I_3$$
thus the Markov transition matrix as well as the escape transition matrix of this concrete interval map are indeed given by $A_f$ and $\widehat{A}_f$ (written above).

The escape dynamics has then an extra symbol $\widehat{2}$ so $\Sigma_{E_f}=\{\widehat{2}\}$, and the alphabet of the escape dynamics is $\{1,2,\hat{2},3,4\}$.
Let $x\in E_f$. Thus $\iota(x)=\widehat{2}$.

For every $V\subset \G_{A_f}^0$, the relative \C-algebra $\Cstar(\G_{A_f},V)$ is generated by 6 partial isometries  and 4 projections. The directed graph $\G_{A_f}$ is as follows:

$$\begin{tikzpicture}[->,>=stealth',shorten >=1pt,auto,node distance=3cm,
                    thick,main node/.style={circle,draw,font=\sffamily\Large\bfseries}]

  \node[main node] (1) {$p_1$};
  \node[main node] (2) [below left of=1] {$p_2$};
  \node[main node] (3) [below right of=2] {$p_4$};
  \node[main node] (4) [below right of=1] {$p_3$};

  \path[every node/.style={font=\sffamily\small}]
    (1) edge node [left] {$s_{13}$} (4)
        edge [right] node[left] {$s_{12}$} (2)
    (2) edge [right] node[left] {$s_{24}$} (3)
    (3) edge [right] node[right] {$s_{43}$} (4)
    (4) edge node {$s_{32}$} (2)
        edge [bend right] node[right] {$s_{31}$} (1);
\end{tikzpicture}
$$
where we label the edges by the defining 6 partial isometries and the vertices by the 4 projections of $\Cstar(\G_{A_f},V)$.

If $x\in E_f$ then we have the concrete 6 partial isometries and 4 projections:
\begin{equation}\label{eqpartprojex1}
  S_{12}, S_{13} S_{24}, S_{31}, S_{32}, S_{43},\ P_1, P_2, P_3, P_4
\end{equation}
given by Eqs.\ (\ref{eqSij}) and (\ref{eqPi}).

If $1\in V$ or $4\in V$, then $\widehat{a}_{i \iota(x)}\not=0$ for some $i\in V$. This implies that partial isometries and projections in \eqref{eqpartprojex1} do not define a representation of $\Cstar(\G_{A_f},V)$.

If $V=\{2,3\}$ the conditions in part (1) of the above Theorem \ref{thmnice100} are fulfilled. In this case Theorem \ref{thmnice100} implies that $\nu_x$ is a faithful representation on $H_x$.

If $V=\{2\}$ then condition (1a) of the above Theorem \ref{thmnice100} is fulfilled. Then the operators in \eqref{eqpartprojex1} generate a (non faithful) representation $\nu_x$ on $H_x$. Similarly with $V=\{3\}$.
\label{example111}
\end{example}

\begin{remark}
Fix $v\in V_2\setminus V_1$ if $V_1\subseteq V_2$. Then we have the relation
$$p_v> \sum_{e\in \G_A^1:\ s(e)=v} s_e s_e^\ast$$
in $\Cstar(\G_A, V_1)$ whereas in $\Cstar(\G_A, V_2)$ we have the equality
$$p_v= \sum_{e\in \G_A^1:\ s(e)=v} s_e s_e^\ast.$$
The quotient map $q: \Cstar(\G_A,V_1)\to \Cstar(\G_A,V_2)$ satisfies
$$q(p_v-\sum_{e\in \G_A^1:\ s(e)=v} s_e s_e^\ast)=0.$$
Thus $q$ is not faithful. In particular, any representation $\pi: \Cstar(\G_A, V_2)\to B(H)$ gives rise to a non-faithful representation $\pi\circ q$ of $\Cstar(\G_A, V_1)$ on the same Hilbert space $H$.
\label{remnoninj}
\end{remark}

Let $A$ be a 0-1 $n\times n$ matrix and $\G_A$ its oriented graph as in \eqref{eqGA}.
We remark that we have a bijection between the vectors $u=(u_1,...,u_n)\in \{0,1\}^n$
and subsets $V$ of $\G_A^0$ such that:
given $u$ let
\begin{equation}\label{equV}
V_u:=\{i: u_i=0\},
\end{equation}
and given $V$, we associate $u_V\in\{0,1\}^n$ with $(u_V)_i=0$ if $i\in V$ and $(u_V)_i=1$ if $i\notin V$. Clearly $V=V_{u_V}$ and $u=u_{V_u}$.

We prove now that indeed we can exhaust all such data $(n; u_1,...,u_n)$ in the representation theory of the relative graph algebras arising from interval maps.
\begin{corollary}
\begin{enumerate}
  \item Let $f\in \mathcal{M}(I)$ and write $\widehat{A}_f$ in the block form
  $$\left[
\begin{array}{c|c}
A_f & B_f \\
\hline
0_{m\times n} & 0_{m\times m}
\end{array}
\right]$$ as in (\ref{eqmatrixa00}). Then each column $\mathbf{u}$ of $B_f$ gives rise to one faithful representation of the relative graph algebra $\Cstar(\G_A,V_{\mathbf{u}})$.

  \item Let $n\in\mathbb{N}$ and $(u_1,...,u_n)\in\{0,1\}^n$. Then there exists an interval map $f\in \mathcal{M}(I)$ and a point $x\in I$
such that 
$$P\widehat{A}_fP^T= \left[
\begin{array}{c|c}
A_f & \begin{array}{c}
 u_1 \\
 \vdots \\
 u_n
\end{array} \\
\hline
0_{1\times n} & 0
\end{array}
\right]$$
for some permutation matrix $P$, and the representation $\nu_x$ of $\Cstar(\G_{A_f},V_u)$ is faithful, where $V_u$ is given as in Eq.\ (\ref{equV}).
\label{corfaithreps}
\end{enumerate}

\end{corollary}
\begin{proof}
Part (1) is a consequence of the definition of $V_{\mathbf{u}}$.
For part (2), we note that the existence of such an interval map $f$ with that escape transition matrix is a consequence of \cite[Proposition 4]{RMP11}. Then take $x\in E_f$. The representation $\nu_x$ of $\Cstar(\G_A,V_u)$ is faithful by Theorem \ref{thmnice100}.
\end{proof}

\begin{remark}
Let $\widehat{A}_f$ be the escape transition matrix of $f$. Take two escape intervals $E_i$ and $E_j$ from $E_f$. Then we have two column vectors $\mathbf{u}_{i}=(u_{i1},...,u_{in})\in \{0,1\}^n$ and   $\mathbf{u}_{j}=(u_{j1},...,u_{jn})\in \{0,1\}^n$ together with the associated sets as in \eqref{equV}: $V_{\mathbf{u}_{i}}$ and $V_{\mathbf{u}_{j}}$.
\begin{enumerate}
  \item Fixing $x\in E_i$ and $y\in E_j$, Theorem \ref{thmnice100} implies that we have two representations $\nu_x$ and $\nu_y$ of the same relative graph algebra $\Cstar(\G_{A_f}, V_{\mathbf{u}_{i}}\cap V_{\mathbf{u}_{j}})$. Note that  $\nu_x$ and $\nu_y$ are faithful representations of  $\Cstar(\G_{A_f}, V_{\mathbf{u}_{i}})$ and  $\Cstar(\G_{A_f}, V_{\mathbf{u}_{j}})$, respectively, by Theorem \ref{thmnice100}.

  \item If $V_{\mathbf{u}_{i}}\not= V_{\mathbf{u}_{j}}$, then $\nu_x$ and $\nu_y$ are not unitary equivalent. Indeed, we then can find $k\in\{1,..., n\}$ such that $\widehat{a}_{k \iota(x)}=1$ and $\widehat{a}_{k \iota(y)}=0$ (or  $\widehat{a}_{k \iota(x)}=0$ and   $\widehat{a}_{k \iota(y)}=1$). Then,
  $$\nu_x (p_k-\sum_{e\in \G_A^1: s(e)=k} s_es_e^\ast)\not=0$$
   whereas
   $$\nu_y (p_k-\sum_{e\in \G_A^1: s(e)=k} s_es_e^\ast)=0$$
    (as in the proof of part (2) in Theorem \ref{thmnice100}). This contradicts $U\nu_x(\cdot) U^\ast =\nu_y(\cdot)$ if they were unitary equivalent. For the irreducibility of $\nu_x$, see \cite{RMP10}.
\end{enumerate}
\end{remark}

If $x\in E_{f}$ then $R_{f}\left(
x\right) $ has a natural structure of a rooted tree. The root of $%
R_{f}\left( x\right) $ is the point $e\left(x\right)$ with no
outgoing edge, so $f^{-1}(e(x))\in \hbox{dom}(f)$ but $e(x)\notin \hbox{dom}(f)$ (see \cite[Remark 3.7(2)]{RMP10}.

\begin{corollary}
Let $\nu_x$ and $\nu_y$ be representations of $\Cstar(\G_{A_f},V)$ as in Definition \ref{defnux}. Then $\nu_x$ and $\nu_y$ are unitary equivalent whenever $R_f(x)$ and $R_f(y)$ are isomorphic as rooted trees. Moreover they are irreducible.
\label{lastcor}
\end{corollary}
\begin{proof}
Since $\nu_x$ is a representation of $\Cstar(\G_{A_f},V)$, $\nu_x\circ q$ is a representation of the Toeplitz algebra $\mathcal{T}_{A_f}$ (which is $\Cstar(\G_{A_f},\emptyset)$), where $q: \mathcal{T}_{A_f}\to \Cstar(\G_{A_f},V)$ is the quotient map. Similar argument goes for $\nu_y\circ q$.

We know by \cite[Theorem 4.5]{RMP10} that the representations $\nu_x\circ q$ and $\nu_y\circ q$ of the Toeplitz algebra $\mathcal{T}_{A_f}$ are unitary equivalent whenever $R_f(x)$ and $R_f(y)$ are isomorphic as rooted trees. If so, then again by definition, the representations $\nu_x$ and $\nu_y$ of $\Cstar(\G_{A_f},V)$ are unitary equivalent.

Finally, By \cite[Theorem 4.7]{RMP10}, $\nu_x$ is irreducible
$\left(\nu_x\circ q(\mathcal{T}_{A_f})\right)^\prime=\mathbb{C}\mathbf{1}$, implying that $\left(\nu_x(\mathcal{T}_{A_f})\right)^\prime=\mathbb{C}\mathbf{1}$ because $q$ is surjective. Therefore $\nu_x$ is irreducible.
\end{proof}

We remark that if $x\in E_f$ and $y\notin E_f$, then by \cite[Proposition 4.6]{RMP10} $\nu_x$ and $\nu_y$ are not unitary equivalent representations of $\Cstar(\G_{A_f},\emptyset)$.

\medskip

{\bf Acknowledgments.}
The authors would like to thank the referees for a close reading of the paper and for useful suggestions.
The first author was partially sponsored by national funds through the Funda\c c\~ao Nacional para a Ci\^encia e Tecnologia, Portugal-FCT, under the project PEst- OE/MAT/UI0117/2014 and CIMA-UE.
The last two authors were partially supported by FCT/Portugal through the project\\ UID/MAT/04459/2013.


\end{document}